\documentclass{article}
\usepackage{amsmath}
\title{The Holonomic Rank of the Fisher-Bingham System of Differential Equations}
\author{Tamio Koyama\footnote{Department of Mathematics, Kobe University
                     and JST CREST Hibi team} ,
        Hiromasa Nakayama${}^*$,  \\
        Kenta Nishiyama\footnote{Graduate School of Information Science
                       and Technology, Osaka University
                       and JST CREST Hibi team} ,
        Nobuki Takayama${}^*$}
\date{February 18, 2013}
\def\comment#1{ }
\def\commentt#1{ }
\def\Del{\partial}
\def\pd#1{ \partial_{#1} }
\def\QED{ Q.E.D. \par \bigbreak} 

\def\C{{\bf C}}
\def\inn{\mathrm{in}}

\newtheorem{theorem}{Theorem}
\newtheorem{proposition}{Proposition}
\newtheorem{lemma}{Lemma}

\begin{document}
\maketitle

{\it Abstract}\/.
The Fisher-Bingham system is a system of linear partial
differential equations satisfied by the Fisher-Bingham 
integral for the $n$-dimensional sphere $S^n$.
The system is given in \cite[Theorem 2]{n3ost2}
and it is shown that it is a holonomic system \cite{koyama1}.
We show that the holonomic rank of the system is equal to $2n+2$.
\bigbreak


\section{Introduction}
Let $x=(x_{ij})$ and $y=(y_i)$ be parameters
such that $x_{ij} = x_{ji}$ for $i\not= j$.
Let $Z$ be a function, which is the normalization constant
of the Fisher-Bingham distribution, defined as 
\begin{equation}  \label{Zdef}
Z(x, y, r) = \int_{S^n(r)} \exp \left(  \sum_{1\leq i \leq j \leq n+1}x_{ij}t_it_j  +\sum_{i=1}^{n+1}y_it_i  \right) |dt|
\end{equation}
where 
$S^n(r) = \{ (t_1, \ldots, t_{n+1}) \mid \sum_{i=1}^{n+1} t_i^2 = r^2, r>0 \}$
is the $n$-dimensional sphere
and $|dt|$ denotes the Haar measure on the sphere.

Let $D$ be the Weyl algebra
$$ D={\bf C}\langle x_{ij}, y_k, r, \pd{ij}, \pd{k}, \pd{r} \mid 
      1 \leq i \leq j \leq n+1, 1 \leq k \leq n+1 \rangle $$
where 
$\pd{ij} = \partial/\partial x_{ij}$,
$\pd{k} = \partial/\partial y_{k}$ and
$\pd{r} = \partial/\partial r$.
It is shown in \cite{koyama1} and \cite{n3ost2} 
that the normalization constant (\ref{Zdef})
of the Fisher-Bingham
distribution is a holonomic function 
in $x, y, r$ and consequently it is annihilated by 
the following holonomic ideal $I$ in $D$
generated by the following operators
in $D$:
\begin{eqnarray*}
\label{diffopA}
&&\Del_{ij} - \Del_i\Del_j, \\ 
\label{diffopB}
&&\sum _{i=1}^{n+1} \Del_i^2 - r^2, \\ 
\label{diffopC}
&&x_{ij}\Del_i^2+2(x_{jj}-x_{ii})\Del _i\Del_j
	-x_{ij}\Del_j^2 \nonumber \\
&& \qquad \qquad
        +\sum _{s \neq i, j}
	\left( x_{sj}\Del_i\Del_s-x_{is}\Del_j\Del_s\right)
	+y_j\Del_i -y_i\Del_j, \\
\label{diffopE}
&& r\Del_r -2\sum _{i\leq j}x_{ij}\Del_i\Del_j
	- \sum_i y_i\Del_i -n.
\label{diffopR}
\end{eqnarray*}
We call the system of differential equations defined by $I$
{\it the Fisher-Bingham system}\/.

For a left ideal $J$ in $D$, 
the holonomic rank of $J$ is defined as the dimension of the 
$K={\bf C}(x_{ij}, y_k, r \mid 1 \leq i \leq j \leq n+1, 1 \leq k \leq n+1)$
vector space 
$K \langle \pd{ij}, \pd{k}, \pd{r} \rangle/K \langle \pd{ij}, \pd{k}, \pd{r} \rangle J$.
The rank is denoted by 
${\rm rank}\,(J)$.
When $J$ is a holonomic ideal, the rank is finite.
The holonomic rank agrees with the dimension of the holomorphic solutions
of the associated system of linear partial differential equations
at generic points and with the size of the Pfaffian equation
associated to $J$.
As to general facts on the holonomic rank, we refer to, e.g., 
the chapters 1 and 2 of \cite{SST}.
The holonomic rank is a fundamental invariant of the $D$-module $D/J$
and there are several attractive studies on holonomic ranks.
For example, Miller, Matusevich and Walther studied holonomic ranks
of $A$-hypergeometric systems
by introducing a new homological method \cite{mmw}.

We are interested in the holonomic rank of the Fisher-Bingham system $I$. 
We prove the following theorem in this paper.
\begin{theorem}  \label{th:main}
$$ {\rm rank}\, (I) = 2n+2 $$
\end{theorem}

In \cite{n3ost2}, we proposed a new method in the statistical inference
which is called the holonomic gradient descent.
The method utilizes a holonomic system of linear partial differential equations
associated to the normalization constant.
The complexity of the method depends on the holonomic rank and
correctness of the method are proved by utilizing the holonomic rank.
In the case of the Fisher-Bingham distribution, 
which is the most fundamental distribution in the directional statistics,
Theorem \ref{th:main} is applied in \cite{knnt1},
which gives a generalization of the result in \cite{n3ost2}
shown with a help of a computer program.
Our method to prove the theorem is the Gr\"obner deformation
to the direction $(-w,w)$,
which is discussed in \cite{SST} for $A$-hypergeometric systems,
and a determination of Gr\"obner bases {\it by hand} with adding
several {\it slack variables}
which do not change the holonomic rank.

\section{The Rank of the Diagonal System} \label{sec:diagonal}

When the matrix $x$ is diagonal, the normalization constant $Z$
satisfies a system of linear partial differential equations for the variables
$x_{ii}$, $y_k$, $r$.
Let $\tilde{I}$
be the left ideal in $D$ generated by
\begin{eqnarray*}
 &&A_i=\Del_{ii} - \Del_{i}^2  \quad  (1\leq i \leq n+1), \\
 &&B    =\sum_{i=1}^{n+1} \Del_i^2 -r^2, \\
 &&C_{ij}=2(x_{ii}-x_{jj})\Del_{i}\Del_{j}+y_i\Del_{j}-y_j\Del_{i}  \quad (1\leq i < j \leq n+1),\\ 
 &&E    =r\Del_r-2\sum_{i=1}^{n+1}x_{ii}\Del_{i}^2-\sum_{i=1}^{n+1}y_i\Del_i-n
\end{eqnarray*}
and 
$ \pd{ij} $, $i\not=j$.
The ideal $\tilde{I}$ annihilates the function $Z$ restricted to the diagonal
of $x$ \cite{knnt1}.

\begin{theorem} \label{th:diagonal_rank}
The holonomic rank of ${\tilde I}$ is $2n+2$.
\end{theorem}

Our proof of the theorem reduces to the proof of the following proposition.
\begin{proposition} \label{pr:fbgr2}
Let R be the ring of differential operators with rational function coefficients
$$R = \C(x_{11}, \ldots, x_{n+1 n+1}, y_1, \ldots, y_{n+1}, r) 
\langle \Del_{11}, \ldots, \Del_{n+1 n+1}, \Del_1, \ldots, \Del_{n+1}, \Del_r \rangle.$$
Let $R \tilde{I}$ be the left ideal of $R$ generated by $A_i, B, C_{ij}, E$.
Let $<$ be the term order on $R$ which is the block order with 
$ \Del_r \gg \{\Del_{ii}\} \gg \{\Del_j\} $.
The order of the block $\{\Del_{ii}\}$ 
is the graded lexicographic order with $\Del_{11} > \cdots > \Del_{n+1 n+1}$ and 
that of the block $\{\Del_{i}\}$ is the graded lexicographic order with 
$\Del_1 > \cdots > \Del_{n+1}$.
A Gr\"{o}bner basis of $R \tilde{I}$ with respect to the term order $<$ is
\begin{align*}
&A_i = \Del_{ii} - \Del_i^2 \quad (i = 1, \ldots, n+1), \\ 
&B = \sum_{i=1}^{n+1} \Del_i^2 - r^2, \\ 
&C_{ij} = 2(x_{ii} - x_{jj}) \Del_i \Del_j + y_i \Del_j - y_j \Del_i \quad (1 \leq i < j \leq n + 1), \\
& (\text{We put } a_{ij} = 2 (x_{ii} - x_{jj}), F_{ij} = y_i \Del_j - y_j \Del_i, 
   C_{ij} = a_{ij} \Del_i \Del_j + F_{ij}), \\
&D_k = \Del_k B - \Del_1 a_{1k}^{-1} C_{1k} - \cdots - \Del_{k-1} a_{k-1 k}^{-1} C_{k-1 k} \quad (k = 1, \ldots, n +1), \\
&E = r \Del_r - 2 \sum_{i=1}^{n+1} x_{ii} \Del_i^2 - \sum_{i=1}^{n+1} y_i \Del_i - n.
\end{align*}
The initial monomials of the Gr\"{o}bner basis are 
\begin{align*}
&  \inn_<(A_i) = \inn_<(\Del_{ii}), 
   \inn_<(B) = \inn_<(\Del_{1})^2, 
   \inn_<(C_{ij}) = \inn_<(\Del_i) \inn_<(\Del_j), \\
&  \inn_<(D_k) = \inn_<(\Del_k)^3, 
   \inn_<(E) = \inn_<(\Del_r).
\end{align*}

\end{proposition}
Here, the initial monomial 
$\inn_<(c(x,y,r)\pd{r}^\gamma \prod \pd{ii}^{\alpha_{ii}} \prod \pd{k}^{\beta_k} )$ 
is defined as the element 
$c(x,y,r) \iota^\gamma \prod \xi_{ii}^{\alpha_{ii}} \prod \eta_{k}^{\beta_k} $
in the polynomial ring with rational function coefficients
$\C(x,y,r)[\xi_{11}, \ldots, \xi_{n+1 n+1}, \eta_1, \ldots, \eta_{n+1}, \iota]$
(see, e.g., \cite[Chapter 1]{SST}).

By the proposition, the standard monomials of the quotient ring
$R/R\tilde{I}$ are
$1, \Del_1, \Del_2, \Del_2^2, \ldots, \Del_{n+1}, \Del_{n+1}^2$.
Therefore, the holonomic rank of $\tilde{I}$ is $2n + 2$ (Theorem
\ref{th:diagonal_rank}).
Let us prove the proposition.

Since $D_k$ is expressed by $B, C_{1k}, \ldots, C_{k-1k}$, 
the operator $D_k$ is the element in $R \tilde{I}$.
In order to prove Proposition \ref{pr:fbgr2}, 
we will show that any $S$-pair for $A_i, B, C_{ij}, D_k, E$ 
is reduced to $0$ by $A_i, B, C_{ij}, D_k, E$. 
The following lemmas are proved by straight forward calculations.

\begin{lemma} \label{lem:nc-cri-d}
Let $P$ and $Q$ be elements in $R$. 
If the initial monomials are coprime, i.e.,  ${\rm gcd}(\inn_<(P), \inn_<(Q)) = 1$, 
then the $S$-pair $S(P, Q)$ is reduced to $[P, Q]$ by $P$ and $Q$
(we denote the reduction by $S(P, Q) \xrightarrow[P, Q]{}^* [P, Q]$), 
where $[P, Q]$ is the commutator of $P$ and $Q$.
In particular, when $[P, Q] = 0$, the $S$-pair $S(P, Q)$ is reduced to $0$.
\end{lemma}

\begin{lemma} \label{lem:cc}
We have 
$$ [A_p, C_{ij}] = 0, $$
$$ [B, C_{ij}] = 0, $$
$$ [C_{ij}, C_{jk}] = C_{ik}, [C_{ij}, C_{ik}] = - C_{jk}, [C_{ik}, C_{jk}] = -C_{ij} \quad (i < j < k), $$ 
$$ [C_{ij}, C_{pq}] = 0 \quad (\{i,j\} \cap \{p,q\} = \emptyset).$$
\end{lemma}

\begin{lemma} \label{lem:cd}
We have 
$$ [D_i, A_j] = 
\begin{cases}
0  & (i < j) \\
2 a_{ji}^{-2} \Del_j C_{ji}  & (i > j), \\ 
\sum_{l=1}^{i-1} 2 a_{li}^{-2} \Del_l C_{li} & (i = j)
\end{cases}$$
$$ [D_i, B] = 0, $$
$$ [D_i, D_j] = -B \Del_j + \sum_{l < i} a_{li}^{-1} a_{ij}^{-1} \Del_l(\Del_iC_{lj} + \Del_l C_{ij}) + \sum_{l < i} a_{li}^{-1} a_{lj}^{-1} (-2\Del_l C_{ij}).$$
When $i, j \neq k$, we obtain 
$$[C_{ij}, D_k] = 
\begin{cases}
0 & (k - 1 < i) \\
0 & (j  \leq k - 1) \\
-a_{ik}^{-1} [C_{ij}, \Del_i C_{ik}] = a_{ik}^{-1} (\Del_i C_{jk} + \Del_j C_{ik}) & (i \leq k - 1 < j)
\end{cases}.$$
\end{lemma}

\begin{lemma} \label{lem:ce}
We have 
\begin{align*}
&[A_i, E] = 0, \\
&[B, E] = -2 B, \\
&[C_{ij}, E] = 0, \\
&[D_i, E] = -3 D_i -2 \sum_{k=1}^{i-1} a_{ki}^{-1} \Del_k C_{ki}.
\end{align*}
\end{lemma}

\noindent
{\it Proof of Proposition \ref{pr:fbgr2}}\/.
We prove that any $S$-pair for $A_i, B, C_{ij}, D_k, E$ is reduced to $0$ 
by $A_i, B, C_{ij}, D_k, E$. 

{\it $S$-pairs of $A_i$ and $A_j, B, C_{ij}$}.
The initial monomials are coprime, and the elements commute.
By Lemma \ref{lem:nc-cri-d}, we obtain
\begin{align*}
& S(A_i, A_j) \longrightarrow^* 0, \\
& S(A_i, B) \longrightarrow^* 0, \\
& S(A_i, C_{jk}) \longrightarrow^* 0. 
\end{align*}

{\it $S$-pair of $B$ and $C_{ij}$}.
When $i > 1$, 
the initial monomials 
$\inn_<(B) = \inn_<(\Del_1)^2$,
$\inn_<(C_{ij}) = \inn_<(\Del_i) \inn_<(\Del_j)$
are coprime.
Operators $B$ and $C_{ij}$ commute by Lemma \ref{lem:cc}.
By Lemma \ref{lem:nc-cri-d},  we have 
$S(B, C_{ij}) \longrightarrow^* 0$.

When $i = 1$, 
the initial monomials 
$\inn_<(B) = \inn_<(\Del_1)^2$,
$\inn_<(C_{1j}) = \inn_<(\Del_1) \inn_<(\Del_j)$
are not coprime.
We obtain the following reduction sequence of the $S$-pair:
\begin{align*}
& S(B, C_{1j}) = a_{1j} \Del_j B - \Del_1 C_{1j} \\ 
& = a_{1j} (\Del_j \Del_2^2 + \cdots + \Del_j^3 + \cdots + \Del_j \Del_{n+1}^2 - \Del_j r^2) - \Del_1 F_{1j} \\
& = a_{1j}((\underline{\Del_j \Del_2^2} + \cdots + \Del_j^3 + \cdots + \Del_j \Del_{n+1}^2 - \Del_j r^2) - \Del_1 a_{1j}^{-1} F_{1j}) \\
& \xrightarrow[C_{2j}]{}^* a_{1j}((\underline{\Del_j \Del_3^2} + \cdots + \Del_j^3 + \cdots + \Del_j \Del_{n+1}^2 - \Del_j r^2) - \Del_1 a_{1j}^{-1} F_{1j} - \Del_2 a_{2j}^{-1} F_{2j}) \\
& \xrightarrow[C_{3j}]{}^* \cdots \xrightarrow[C_{j-1j}]{}^* a_{1j} D_j \xrightarrow[D_j]{} 0.
\end{align*}

{\it $S$-pair of $C_{ij}$ and $C_{kl}$ $(\{i,j\} \cap \{k,l\} = \emptyset)$}.
The initial monomials 
$\inn_<(C_{ij}) = \inn_<(\Del_i) \inn_<(\Del_j), \inn_<(C_{kl}) = \inn_<(\Del_k) \inn_<(\Del_l)$ 
are coprime.
Operators $C_{ij}$ and $C_{kl}$ commute by Lemma \ref{lem:cc}.
By Lemma \ref{lem:nc-cri-d}, we obtain
$$ S(C_{ij}, C_{kl}) \xrightarrow[C_{ij}, C_{kl}]{}^* 0.$$

{\it $S$-pair of $C_{ij}$ and  $C_{jk}$ $(i < j < k)$}.
We have the following reduction sequence of the $S$-pair: 
\begin{align*}
S(C_{ij}, C_{jk}) &= a_{jk} \Del_k C_{ij} - a_{ij} \Del_i C_{jk} = -a_{ik} y_j \Del_i \Del_k + a_{jk} y_i \Del_j \Del_k + a_{ij} y_k \Del_i \Del_j \\
&\xrightarrow[C_{ij}, C_{jk}, C_{kl}]{}^* y_i(-F_{jk}) + y_k (-F_{ij}) - y_j(-F_{ik}) = 0.
\end{align*}
The $S$-pair of $C_{ij}$ and $C_{ik}$ and that of $C_{ik}$ and $C_{jk}$ are also reduced to $0$.

{\it $S$-pair of $D_i$ and $A_j$}.
The initial monomials 
$\inn_<(D_i) = \inn_<(\Del_i)^3, \inn_<(A_j) = \inn_<(\Del_{jj})$ are coprime.
When $i < j$, 
operators $D_i$ and $A_j$ commute. By Lemma \ref{lem:nc-cri-d},  we have
$$ S(D_i, A_j) \xrightarrow[D_i, A_j]{}^* 0.$$

When $i > j$,  by Lemmas \ref{lem:nc-cri-d} and \ref{lem:cd}, we have
$$ S(D_i, A_j) \xrightarrow[D_i, A_j]{}^* [D_i, A_j] = 2 a_{ji}^{-2} \Del_j C_{ji} \xrightarrow[C_{ji}]{}^* 0.$$

When $i = j$, by Lemmas \ref{lem:nc-cri-d} and \ref{lem:cd}, we have 
$$ S(D_i, A_i) \xrightarrow[D_i, A_j]{}^* [D_i, A_i] = \sum_{l=1}^{i-1} 2 a_{li}^{-2} \Del_l C_{li} \xrightarrow[C_{1i}, \ldots, C_{i-1 i}]{}^* 0.$$ 

{\it $S$-pair of $D_i$ and  $B$}.
When $i = 1$, the $S$-pair is $S(D_1, B) = \Del_1 B - \Del_1 B = 0$.

When $i > 1$, 
the initial monomials $\inn_<(D_i) = \inn_<(\Del_i)^3$,
$\inn_<(B) = \inn_<(\Del_1)^2$ are coprime. 
By Lemmas \ref{lem:nc-cri-d} and  \ref{lem:cd},   we have
$$ S(D_i, B) \xrightarrow[D_i, B]{}^* [D_i, B] = 0.$$

{\it $S$-pair of $D_i$ and $D_j$}.
The initial monomials $\inn_<(D_i) = \inn_<(\Del_i)^3$,
$\inn_<(D_j) = \inn_<(\Del_j)^3$ are coprime.
By Lemmas \ref{lem:nc-cri-d} and \ref{lem:cd},  we have
$$ S(D_i, D_j) \xrightarrow[D_i, D_j]{}^* [D_i, D_j] \xrightarrow[B, C_{1j}, \ldots, C_{ij}]{}^* 0.$$

{\it $S$-pair of $C_{ij}$ and $D_k$}.
The initial monomials are $\inn_<(C_{ij}) = \inn_<(\Del_i)\inn_<(\Del_j)$,
$\inn_<(D_k) = \inn_<(\Del_k)^3$.
When $i \neq k$ and $j \neq k$, 
the initial monomials are coprime.
By Lemmas \ref{lem:nc-cri-d} and  \ref{lem:cd}, we have 
$$ S(C_{ij}, D_k) \xrightarrow[C_{ij}, D_k]{}^* [C_{ij}, D_k] \xrightarrow[C_{ik}, C_{jk}]{}^* 0.$$

When $i = k$, 
the initial monomials $\inn_<(C_{ij}) = \inn_<(\Del_i) \inn_<(\Del_j), \inn_<(D_i) = \inn_<(\Del_i)^3$ 
are not coprime.
This case needs a care of an order of applying reductions.
We will reduce the $S$-pair by $D_j$ and then reduce remainders by $C_{ij}$'s. 
\begin{align*}
S(C_{ij}, D_i) &= \Del_i^2 C_{ij} - a_{ij} \Del_j D_i  \\
&= \Del_i^2 F_{ij} - a_{ij} \underline{\Del_i \Del_j^3} - a_{ij} \Del_j (\sum_{l=i+1, l \neq j}^{n+1} \Del_i \Del_{l}^2 - \Del_i r^2 - \sum_{l=1}^{i-1} \Del_l a_{li}^{-1} F_{li})  \\
&\xrightarrow[D_j]{}^* \Del_i^2 F_{ij} + a_{ij}
 \Del_i(\sum_{l=j+1}^{n+1} \Del_j \Del_l^2 - \Del_j r^2 -
 \sum_{l=1}^{j-1} \Del_l a_{lj}^{-1} F_{lj}) \\
&\qquad\qquad\qquad\qquad -a_{ij} \Del_j (\sum_{l=i+1, l \neq j}^{n+1} \Del_i \Del_{l}^2 - \Del_i r^2 - \sum_{l=1}^{i-1} \Del_l a_{li}^{-1} F_{li})  \\
&= - a_{ij} \Del_i (\Del_{i+1} a_{i+1j}^{-1} C_{i+1j} + \cdots + \Del_{j-1} a_{j-1j}^{-1} C_{j-1j}) \\
&\qquad\qquad\qquad\qquad - a_{ij} \Del_i (\sum_{l=1}^{i-1} a_{lj}^{-1} F_{lj}) + a_{ij} \Del_j (\sum_{l=1}^{i-1} \Del_l a_{li}^{-1} F_{li}) \\
&\xrightarrow[C_{i+1j}, \ldots, C_{j-1j}]{}^* - a_{ij} \Del_i (\sum_{l=1}^{i-1} a_{lj}^{-1} F_{lj}) + a_{ij} \Del_j (\sum_{l=1}^{i-1} \Del_l a_{li}^{-1} F_{li}) \\
&= -a_{ij} \sum_{l=1}^{i-1} (a_{lj}^{-1} \Del_l \Del_i F_{lj} - a_{li}^{-1} \Del_l \Del_j F_{li}). 
\end{align*}
Since 
$a_{lj}^{-1} \Del_l \Del_i F_{lj} - a_{li}^{-1} \Del_l \Del_j F_{li} 
\xrightarrow[C_{lj}, C_{li}, C_{ij}]{}^* 0$, 
the $S$-pair 
$S(C_{ij}, D_i)$ is reduced to $0$. 

When $j = k$, 
the initial monomials $\inn_<(C_{ij}) = \inn_<(\Del_i) \inn_<(\Del_j), \inn_<(D_j) = \inn_<(\Del_j)^3$ 
are not coprime.
This case also needs a care of an order of applying reductions.
\begin{align*}
S(C_{ij}, D_j) &= \Del_j^2 C_{ij} - a_{ij} \Del_i D_j  \\
&= F_{ij} (\Del_i^2  + \Del_j^2) - a_{ij} \Del_i 
   (\sum_{l=j+1}^{n+1} \Del_j \Del_l^2 - \Del_j r^2 - \sum_{l=1, l \neq i}^{j-1} \Del_l a_{lj}^{-1} F_{lj})  \\
&\xrightarrow[B]{}^* F_{ij} (-\sum_{l=1, l \neq i, j}^{n+1} \Del_l^2 + r^2) - a_{ij} \Del_i 
   (\sum_{l=j+1}^{n+1} \Del_j \Del_l^2 - \Del_j r^2 - \sum_{l=1, l \neq i}^{j-1} \Del_l a_{lj}^{-1} F_{lj})  \\
&= F_{ij} (-\sum_{l=1, l \neq i, j}^{n+1} \Del_l^2) + r^2 (a_{ij} \Del_i \Del_j + F_{ij})  - a_{ij} \Del_i 
   (\sum_{l=j+1}^{n+1} \Del_j \Del_l^2 - \sum_{l=1, l \neq i}^{j-1} \Del_l a_{lj}^{-1} F_{lj})  \\
&\xrightarrow[C_{ij}]{}^* F_{ij} (-\sum_{l=1, l \neq i, j}^{n+1} \Del_l^2) - a_{ij} \Del_i 
   (\sum_{l=j+1}^{n+1} \Del_j \Del_l^2 - \sum_{l=1, l \neq i}^{j-1} \Del_l a_{lj}^{-1} F_{lj}) \\
&= (\sum_{l=j+1}^{n+1} \Del_l^2) (-a_{ij} \Del_i \Del_j - F_{ij}) - F_{ij} \sum_{l=1, l \neq i}^{j-1} \Del_l^2 + 
    a_{ij} \Del_i \sum_{l=1, l \neq i}^{j-1} \Del_l a_{lj}^{-1} F_{lj}  \\  
&\xrightarrow[C_{ij}]{}^* - F_{ij} \sum_{l=1, l \neq i}^{j-1} \Del_l^2 + 
    a_{ij} \Del_i \sum_{l=1, l \neq i}^{j-1} \Del_l a_{lj}^{-1} F_{lj} \\
&= \sum_{l=1, l \neq i}^{j-1} \Del_l a_{ij} (-a_{ij}^{-1} \Del_l F_{ij} + a_{lj}^{-1} \Del_i F_{lj}).
\end{align*}
Since 
$-a_{ij}^{-1} \Del_l F_{ij} + a_{lj}^{-1} \Del_i F_{lj}
\xrightarrow[C_{lj}, C_{li}, C_{ij}]{}^* 0$, 
the $S$-pair 
$S(C_{ij}, D_j)$ is reduced to $0$. 

{\it $S$-pair of $E$ and $A_i$}.
The initial monomials $\inn_<(E) = \inn_<(\Del_r), \inn_<(A_i) = \inn_<(\Del_{ii})$ are coprime.
By Lemmas \ref{lem:nc-cri-d} and \ref{lem:ce}, we have 
$$ S(A_i, E) \xrightarrow[A_i, E]{}^* [A_i, E] = 0.$$

{\it $S$-pair of $E$ and $B$}.
The initial monomials 
$\inn_<(E) = \inn_<(\Del_r), \inn_<(B) = \inn_<(\Del_{1})^2$ 
are coprime.
By Lemmas \ref{lem:nc-cri-d} and \ref{lem:ce}, we have
$$ S(B, E) \xrightarrow[B, E]{}^* [B, E] = -2B \xrightarrow[B]{}^* 0.$$

{\it $S$-pair of $E$ and $C_{ij}$}.
The initial monomials 
$\inn_<(E) = \inn_<(\Del_r), \inn_<(C_{ij}) = \inn_<(\Del_i) \inn_<(\Del_j)$ 
are coprime.
By Lemmas \ref{lem:nc-cri-d} and \ref{lem:ce}, we have
$$ S(C_{ij}, E) \xrightarrow[C_{ij}, E]{}^* [C_{ij}, E] = 0.$$

{\it $S$-pair of $D_i$ and $E$}.
The initial monomials 
$\inn_<(E) = \inn_<(\Del_r), \inn_<(D_i) = \inn_<(\Del_i)^3$ 
are coprime.
By Lemmas \ref{lem:nc-cri-d} and \ref{lem:ce}, we have
$$ S(D_i, E) \xrightarrow[D_i,E]{}^* [D_i, E] 
   = -3 D_i -2 \sum_{k=1}^{i-1} a_{ki}^{-1} \Del_k C_{ki}
      \xrightarrow[D_i, C_{1i}, \ldots, C_{i-1i}]{}^* 0.$$

We have proved that any $S$-pair is reduced to $0$.
By Buchberger's criterion, 
the set $\{A_i, B, C_{ij}, D_k, E\}$ is a Gr\"{o}bner basis of $R\tilde{I}$.
\QED

\section{Gr\"obner Deformation of the Fisher-Bingham System}
\label{sec:deform}
Consider the system of differential equations $I \cdot f = 0$.
Intuitively speaking, we want to prove that the system $I$ can be
deformed to the diagonal system ${\tilde I}$ without
increasing the holonomic rank.
This can be done by a Gr\"obner basis computation with a weight
vector $(-w,w)$ \cite[Theorem 2.2.1]{SST}.
However a straightforward calculation does not seem to be easy.
We need to use some technical tricks to determine a suitable
Gr\"obner deformation.
Since these tricks may look too technical for the general $n$,
we explain them in the case of $n=1$ in Section \ref{sec:n=1}
to clarify our idea without technical details of this section.
Readers are expected to refer to the Section \ref{sec:n=1}
when technicalities get complicated.
 
We will introduce new indeterminates to make our Gr\"obner basis computation
possible by hand
with employing the idea of the proof of \cite[Theorem 3.1.3]{SST}.
Let $a_{pq}$, $b_i$, $c_i$, $d$
($1 \leq p \leq q \leq n+1$, $1 \leq i \leq n+1$) be constants,
which we call slack variables when they are regarded as indeterminates.
We put $g = \left (r^{d^3}\prod_{p\leq q} x_{pq}^{a_{pq}^3} \right)  f$
and make a change of the independent variables $y_i$ by $y_i+b_i c_i$.
Then, the system of differential equations for the function $g$
is 
$I' \cdot g = 0$
where $I'$  is the left ideal in $D$ generated by the set of operators
$G'=\{A_{pq}', B, C_{ij}', E'\}$ where
\begin{align*}
&A_{pq}' = x_{pq}\Del_{pq} - x_{pq}\Del_p\Del_q -a_{pq}^3, \\ 
&B = \sum _{i=1}^{n+1} \Del_i^2 - r^2, \\ 
&C_{ij}' = x_{ij}\Del_i^2+2(x_{jj}-x_{ii})\Del _i\Del_j
	-x_{ij}\Del_j^2 \nonumber \\
&\qquad\quad
        +\sum _{s \neq i, j}
	\left( x_{sj}\Del_i\Del_s-x_{is}\Del_j\Del_s\right)
	+(y_j+b_j c_j) \Del_i - (y_i+b_i c_i) \Del_j, \\
&E' = r\Del_r -2\sum _{i\leq j}x_{ij}\Del_i\Del_j
	- \sum_{i=1}^{n+1} (y_i+b_i c_i) \Del_i -n -d^3.
\end{align*}
The key fact is that
the holonomic rank of $I$ for $f$ agrees with
the holonomic rank of $I'$ for $g$
for any constants $a_{pq}$, $b_i$, $c_i$, $d$.

Let us make the same change of the variables
for the diagonal system;
let ${\tilde I}'$ be the left ideal of $D$ generated by the set of operators
$\tilde{G}'=\{\tilde{A}_{ii}', B, \tilde{C}_{ij}', \tilde{E}',
x_{ij}\pd{ij}-a_{ij}^3 \: (i\not=j)\}$ where
\begin{align*}
 &\tilde{A}_{ii}' = x_{ii}\Del_{i}^2-x_{ii}\Del_{ii}-a_{ii}^3  \quad  (1\leq i \leq n+1), \\
 &B = \sum_{i=1}^{n+1}\Del_{i}^2-r^2, \\
 &\tilde{C}_{ij}' = 2(x_{jj}-x_{ii})\Del_{i}\Del_{j}
   +(y_j+b_j c_j)\Del_{i}-(y_i+b_ic_i)\Del_{j} \quad (1\leq i < j \leq n+1),\\ 
 &\tilde{E}' = r\Del_r-2\sum_{i=1}^{n+1}x_{ii}\Del_{i}^2
    -\sum_{i=1}^{n+1}(y_i+b_i c_i) \Del_i-n-d^3.
\end{align*}
The holonomic ranks of ${\tilde I}$ and ${\tilde I}'$ agree.

Define the weight vector $w$
by $w_{ij} = 1$, $(i \not= j)$,
$w_{ii} = 0$, $w_k =0$, and $w_r=0$.
Here, $w_{ij}$ stands for $\pd{ij}$,
$w_k$ stands for $\pd{k}$, and
$w_r$ stands for $\pd{r}$.
The initial form ${\rm in}_{(-w,w)}(\ell)$, $\ell \in D$,
is the sum of the highest $(-w,w)$-degree terms in $\ell$
and 
${\rm in}_{(-w,w)}(I')$ is the left ideal
generated by $\ell \in I'$
where the weight $-w$ stands for space variables 
$x_{ij}$, $y_k$, $r$
corresponding to differential operators
$\pd{ij}$, $\pd{k}$, $\pd{r}$ respectively
(see, e.g., \cite[Chapter 1]{SST}).

\begin{theorem}  \label{th:grdeformation}
For generic complex numbers $a_{pq}$, $b_i$, $c_i$, $d$, we have
\begin{equation}  \label{eq:grdeformation}
{\rm in}_{(-w,w)} (I') =  {\tilde I}'.
\end{equation}
\end{theorem}
In order to prove the theorem,
we regard $a_{pq}$, $b_i$, $c_i$, $d$ as ring variables with the weight 0 and
consider the following homogenized system of $I'$:
\begin{align*}
&A_{pq}'^h = hx_{pq}\Del_{pq} - x_{pq}\Del_p\Del_q -a_{pq}^3, \\ 
&B = \sum _{i=1}^{n+1} \Del_i^2 - r^2, \\ 
&C_{ij}'^h = x_{ij}\Del_i^2+2(x_{jj}-x_{ii})\Del _i\Del_j
	-x_{ij}\Del_j^2 \nonumber \\
&\qquad\quad
        +\sum _{s \neq i, j}
	\left( x_{sj}\Del_i\Del_s-x_{is}\Del_j\Del_s\right)
	+(hy_j+b_j c_j) \Del_i - (hy_i+b_i c_i) \Del_j, \\
&E'^h = hr\Del_r -2\sum _{i\leq j}x_{ij}\Del_i\Del_j
	- \sum_{i=1}^{n+1} (hy_i+b_i c_i) \Del_i -nh^3 -d^3.
\end{align*}
Let $I'^h$ be a left $D^h[a,b,c,d]$-ideal generated by the set of operators
$G'^h:=\{A_{pq}'^h, B, C_{ij}'^h, E'^h\}$,
where $D^h[a,b,c,d]$ is the homogenized Weyl algebra of the ring
$D[a,b,c,d]={\bf C}[a,b,c,d]
\langle x_{ij}, y_k, r, \pd{ij}, \pd{k}, \pd{r} \rangle$
with the homogenization variable $h$
(see, e.g., \cite[Section 9]{ot}).
We introduce a new term order $<_{(-w,w,0)}^h$ over $D^h[a,b,c,d]$,
which compares the total degree first, $(-w,w,0)$-degree second,
otherwise we apply the following block order as a tie breaker:
$d \gg r \gg \{a_{pq} \mid i\leq j\} \gg \{b_k\} \gg \{c_k\} \gg \{y_{k}\}
\gg \pd{r} \gg \{\pd{ij} \mid i<j\} \gg \{\pd{ii}\} \gg \{\pd{k}\}
\gg \{x_{ij} \mid i<j\} \gg \{x_{ii}\} \gg h$.
Here, the block $\{b_k\}$ has a lexicographic order so that
$b_1 > b_2> \cdots > b_{n+1}$.
The use of this tie breaker is a key of our calculation.
Although, the initial monomial ${\rm in}_{<_{(-w,w,0)}^h}(\ell)$
is an element of the associated commutative ring,
we denote it by the associated element in $D^h[a,b,c,d]$
as long as no confusion arises.
For example, we denote ${\rm in}_{<_{(-w,w,0)}^h}(\pd{ij})$
by $\pd{ij}$ instead of $\xi_{ij}$.

\begin{proposition} \label{prop:homogGb}
A Gr\"obner basis of $I'^h$ with respect to $<_{(-w,w,0)}^h$ is 
$G'^h=\{ A_{pq}'^h,B,C_{ij}'^h,E'^h \}$.
\end{proposition}

We need four lemmas for proving the proposition.
These can be obtained by a straightforward calculation.

\begin{lemma} \label{lem:ini-h}
The initial monomials of the generators of $I'^h$ are
\begin{enumerate}
\item ${\rm in}_{<_{(-w,w,0)}^h}(A_{pq}'^h) = -a_{pq}^3$,
\item ${\rm in}_{<_{(-w,w,0)}^h}(B) = -r^2$,
\item ${\rm in}_{<_{(-w,w,0)}^h}(C_{ij}'^h) = -b_ic_i\pd{j}$,
\item ${\rm in}_{<_{(-w,w0)}^h}(E'^h) = -d^3$.
\end{enumerate}
In particular, they are pairwise coprime expect in the case that
the pair
${\rm in}_{<_{(-w,w,0)}^h}(C_{ij}'^h)$ and ${\rm in}_{<_{(-w,w,0)}^h}(C_{ik}'^h)$
and the pair
${\rm in}_{<_{(-w,w,0)}^h}(C_{ik}'^h)$ and ${\rm in}_{<_{(-w,w,0)}^h}(C_{jk}'^h)$.
\end{lemma}

\begin{lemma} \label{lem:com-h}
The commutators of two generators of $I'^h$ are
\begin{enumerate}
\item $[A_{pq}'^h,A_{rs}'^h] = [A_{pq}'^h,B] = [A_{pq}'^h,C_{ij}'^h]
       = [A_{pq}'^h,E'^h] = 0$ \: $(\text{for any } p,q,r,s,i,j)$,
\item $[B,C_{ij}'^h] = 0$, $[B,E'^h] = -2hB$,
\item $[C_{ij}'^h,C_{kl}'^h] = 0$ \: $(\{i,j\} \cap \{k,l\} = \emptyset)$,
\item $[C_{ij}'^h,C_{jk}'^h] = C_{ki}'^h := -C_{ik}'^h$,
$[C_{ij}'^h,C_{ik}'^h] = C_{jk}'^h$,
$[C_{ik}'^h,C_{jk}'^h] = C_{ij}'^h$ \: $(i<j<k)$,
\item $[C_{ij}'^h,E'^h]=0$.
\end{enumerate}
\end{lemma}

\begin{lemma} \label{lem:s-pair-h}
The following holds for $S$-pairs of generators of $I'^h$:
\begin{enumerate}
\item $S(A_{pq}'^h,A_{rs}'^h), S(A_{pq}'^h,B), S(A_{pq}'^h,C_{ij}'^h),
       S(A_{pq}'^h,E'^h) \longrightarrow^* 0$ \: $(\text{for any } p,q,r,s,i,j)$,
\item $S(B,C_{ij}'^h) \longrightarrow^* 0$, $S(B,E'^h) \longrightarrow^* 0$,
\item $S(C_{ij}'^h,C_{kl}'^h) \longrightarrow^* 0$ \: $(\{i,j\} \cap \{k,l\} = \emptyset)$,
\item $S(C_{ij}'^h,C_{jk}'^h) \longrightarrow^* 0$ \: $(i<j<k)$,
\item $S(C_{ij}'^h,E'^h) \longrightarrow^* 0$.
\end{enumerate}
\end{lemma}

\noindent
{\it Proof}\/.
Lemma \ref{lem:nc-cri-d} holds in the homogenized Weyl algebra
$D^h[a,b,c,d]$ too.
Therefore, these follow from Lemmas \ref{lem:ini-h} and \ref{lem:com-h}.
\QED

\begin{lemma} \label{lem:cycle-h}
We put $\hat{C}_{ij}'^h := b_jc_j\pd{i}-b_ic_i\pd{j}$ and
$\check{C}_{ij}'^h := C_{ij}'^h-\hat{C}_{ij}'^h$.
Then the following cyclic relations hold:
\begin{enumerate}
\item $\pd{k}\hat{C}_{ij}'^h+\pd{i}\hat{C}_{jk}'^h+\pd{j}\hat{C}'^h_{ki}=\hat{C}_{ij}'^h\pd{k}+\hat{C}_{jk}'^h\pd{i}+\hat{C}_{ki}'^h\pd{j}=0$, \label{lem:cycle-h:1}
\item $\pd{k}\check{C}_{ij}'^h+\pd{i}\check{C}_{jk}'^h+\pd{j}\check{C}_{ki}'^h=\check{C}_{ij}'^h\pd{k}+\check{C}_{jk}'^h\pd{i}+\check{C}_{ki}'^h\pd{j}=0$, \label{lem:cycle-h:2}
\item $\pd{k}C_{ij}'^h+\pd{i}C_{jk}'^h+\pd{j}C_{ki}'^h=C_{ij}'^h\pd{k}+C_{jk}'^h\pd{i}+C_{ki}'^h\pd{j}=0$, \label{lem:cycle-h:3}
\item $b_kc_k\hat{C}_{ij}'^h+b_ic_i\hat{C}_{jk}'^h+b_jc_j\hat{C}_{ki}'^h=\hat{C}_{ij}'^hb_kc_k+\hat{C}_{jk}'^hb_ic_i+\hat{C}_{ki}'^hb_jc_j=0$. \label{lem:cycle-h:4}
\end{enumerate}
\end{lemma}

\noindent
{\it Proof of Proposition \ref{prop:homogGb}}\/.
By Lemma \ref{lem:s-pair-h}, we only need to check that the $S$-pairs
$S(C_{ij}'^h,C_{ik}'^h)$ and $S(C_{ik}'^h,C_{jk}'^h)$ are reduced to zero.

The former is
$S(C_{ij}'^h,C_{ik}'^h) = \pd{k}C_{ij}'^h-\pd{j}C_{ik}'^h
= \pd{k}C_{ij}'^h+\pd{j}C_{ki}'^h$, 
because the initial monomials are
${\rm in}_{<_{(-w,w,0)}^h}(C_{ij}'^h) = -b_ic_i\pd{j}$ and
${\rm in}_{<_{(-w,w,0)}^h}(C_{ik}'^h) = -b_ic_i\pd{k}$.
The $S$-pair is equal to $-\pd{i}C_{jk}'^h$ by the \ref{lem:cycle-h:3}rd formula of Lemma \ref{lem:cycle-h}.
This implies that it is reduced to zero.

The latter is
$S(C_{ik}'^h,C_{jk}'^h) = b_jc_jC_{ik}'^h-b_ic_iC_{jk}'^h$,
because the initial monomials are
${\rm in}_{<_{(-w,w,0)}^h}(C_{ik}'^h) = -b_ic_i\pd{k}$ and
${\rm in}_{<_{(-w,w,0)}^h}(C_{jk}'^h) = -b_jc_j\pd{k}$.

Firstly, we show that it can be expressed as
\begin{align} \label{eq:std-rep}
S(C_{ik}'^h,C_{jk}'^h) &=
  \left (\sum_{s=1}^{n+1} (\delta_{ks}+1)x_{ks}\pd{s} + hy_k+b_kc_k \right ) C_{ij}'^h \nonumber \\
&+\left (\sum_{s=1}^{n+1} (\delta_{is}+1)x_{is}\pd{s} + hy_i \right ) C_{jk}'^h \nonumber \\
&+\left (\sum_{s=1}^{n+1} (\delta_{js}+1)x_{js}\pd{s} + hy_j \right ) C_{ki}'^h.
\end{align}
Here, the symbol $\delta$ is the Kronecker delta.
This expression of the $S$-pair is a key of our proof.
Since each monomial of the left hand side (LHS in short)
has at least one variable in $b_i$, $b_j$ or $b_k$,
we divide the right hand side (RHS in short) into two parts:
$S_1$ contains these variables, $S_2$ does not contain them.
Then,
\begin{align*}
 S_2
&= \left (\sum_{s=1}^{n+1} (\delta_{ks}+1)x_{ks}\pd{s} + hy_k \right ) \check{C}_{ij}'^h
  +\left (\sum_{s=1}^{n+1} (\delta_{is}+1)x_{is}\pd{s} + hy_i \right ) \check{C}_{jk}'^h \\
&+ \left (\sum_{s=1}^{n+1} (\delta_{js}+1)x_{js}\pd{s} + hy_j \right ) \check{C}_{ki}'^h
\end{align*}
is equal to zero by a straightforward calculation and we have
\begin{align*}
S_1 &=
b_kc_kC_{ij}'^h + \left (\sum_{s=1}^{n+1} (\delta_{ks}+1)x_{ks}\pd{s} + hy_k \right ) \hat{C}_{ij}'^h\\
&\qquad +\left (\sum_{s=1}^{n+1} (\delta_{is}+1)x_{is}\pd{s} + hy_i \right ) \hat{C}_{jk}'^h
+\left (\sum_{s=1}^{n+1} (\delta_{js}+1)x_{js}\pd{s} + hy_j \right ) \hat{C}_{ki}'^h\\
&= b_kc_kC_{ij}'^h\\
&+b_jc_j \left (\sum_{s=1}^{n+1} (\delta_{ks}+1)x_{ks}\pd{s}\pd{i} + hy_k\pd{i} \right )
 -b_ic_i \left (\sum_{s=1}^{n+1} (\delta_{ks}+1)x_{ks}\pd{s}\pd{j} + hy_k\pd{j} \right )\\
&+b_kc_k \left (\sum_{s=1}^{n+1} (\delta_{is}+1)x_{is}\pd{s}\pd{j} + hy_i\pd{j} \right )
 -b_jc_j \left (\sum_{s=1}^{n+1} (\delta_{is}+1)x_{is}\pd{s}\pd{k} + hy_i\pd{k} \right )\\
&+b_ic_i \left (\sum_{s=1}^{n+1} (\delta_{js}+1)x_{js}\pd{s}\pd{k} + hy_j\pd{k} \right )
 -b_kc_k \left (\sum_{s=1}^{n+1} (\delta_{js}+1)x_{js}\pd{s}\pd{i} + hy_j\pd{i} \right )\\
&=b_jc_jC_{ik}'^h-b_ic_iC_{jk}'^h = S(C_{ik}'^h,C_{jk}'^h).
\end{align*}

Finally, we show that the expression (\ref{eq:std-rep}) is a standard representation.
We note that the RHS of (\ref{eq:std-rep}) has three terms,
which appear in the first, in the second and in the third lines 
of (\ref{eq:std-rep}) respectively.
We may show that the initial monomial of the LHS is no less than
the initial monomials of the three terms of the RHS with respect to the term order $<_{(-w,w,0)}^h$.
The initial monomial of the LHS $b_ic_ib_kc_k\pd{j}$ coincides with
that of the 1st term of the RHS.
Moreover, all monomials appearing in the 2nd and the 3rd term of the RHS have
the same total degree 5 and $(-w,w,0)$-degree 0, and they have a degree
at most one with respect to $b_i$, $b_j$ and $b_k$.
This implies that they are less than $b_ic_ib_kc_k\pd{j}$.
Thus, we have shown that the expression is a standard representation
and the $S$-pair is reduced to zero.

Since all $S$-pairs are reduced to zero,
the conclusion follows from Buchberger's criterion. \QED

\noindent
{\it Proof of Theorem \ref{th:grdeformation}}\/.
Let $I'(a,b,c,d)$ be the left ideal generated by $G'$
in the ring $D[a,b,c,d]$.
It follows from Proposition \ref{prop:homogGb} and \cite[Theorem 1.2.4]{SST}
that $G'^h|_{h=1}=G'$ is a Gr\"obner basis of $I'(a,b,c,d)$
with respect to $<_{(-w,w,0)}$.
Therefore, we have
$${\rm in}_{<_{(-w,w,0)}}(I'(a,b,c,d))
= \langle {\rm in}_{<_{(-w,w,0)}}(G') \rangle
= \langle \tilde{G}' \rangle \quad \text{ in } D[a,b,c,d].$$
We may replace $D[a,b,c,d]$ by
$D(a,b,c,d) = {\bf C}(a,b,c,d) \langle x_{ij}, y_k, r, \pd{ij}, \pd{k}, \pd{r} \rangle$.
Here, ${\bf C}(a,b,c,d)$ is the rational function field
with variables $a=(a_{pq}), b=(b_i), c=(c_i)$ and $d$.
Then, the following holds:
$${\rm in}_{(-w,w)}(I'(a,b,c,d)) =
\langle \tilde{G}' \rangle \quad \text{ in } D(a,b,c,d).$$
This implies the conclusion.
\QED

\noindent
{\it Proof of the Main theorem \ref{th:main}}\/.
It follows from (\ref{eq:grdeformation}) that
${\rm rank}\,(I) \geq {\rm rank}\, ({\rm in}_{(-w,w)}(I'))$
by Theorem 2.2.1 in \cite{SST}.
Therefore, we have ${\rm rank}\,(I) \geq 2n+2$
by Theorem \ref{th:diagonal_rank}.
The opposite inequality follows from Theorem 3 of \cite{n3ost2}.
\QED

\section{A Proof in the case of $n=1$}  \label{sec:n=1}
In order to clarify ideas of the proof in Section \ref{sec:deform},
we present a proof of our theorem in the case of $n=1$.

The $1$-dimensional Fisher-Bingham system of differential equations
$I \subset \C\langle x_{11},x_{12},x_{22},y_1,y_2,r,
\pd{11},\pd{12},\pd{22},\pd{1},\pd{2},\pd{r} \rangle$
is
\begin{align*}
I =\langle &A_{11}=\pd{11}-\pd{1}^2, \; A_{12}=\pd{12}-\pd{1}\pd{2},\\
& A_{22}=\pd{22}-\pd{2}^2, \\
& B=\pd{1}^2+\pd{2}^2-r^2, \\  
& C_{12}=x_{12}\pd{1}^2+2(x_{22}-x_{11})\pd{1}\pd{2}-x_{12}\pd{2}^2
  +y_2\pd{1}-y_1\pd{2}, \\
& E=r\pd{r}-2(x_{11}\pd{1}^2+x_{12}\pd{1}\pd{2}+x_{22}\pd{2}^2)
-(y_1\pd{1}+y_2\pd{2})-1 \rangle.
\end{align*}
The upper bound of ${\rm rank}(I)$ is $2\cdot1+2 = 4$
as given in \cite[Theorem 3]{n3ost2}.
We will show the lower bound of the
${\rm rank}(I)$ is $4$ by using the following general inequality
\cite[Theorem 2.2.1]{SST}:
$${\rm rank}(I) \geq {\rm in}_{(-w,w)}(I).$$

Firstly, we make some change of variables,
because it seems to be difficult to calculate ${\rm in}_{(-w,w)}(I)$ directly.
Let $a_{11},a_{12},a_{22},b_1,b_2,c_1,c_2,d$ be constants,
which will be used as slack variables.
We put
$g=r^{d^3}x_{11}^{a_{11}^3}x_{12}^{a_{12}^3}x_{22}^{a_{22}^3}f$
where the function $f$ is a solution of the system of
differential equations $I\cdot f=0$.
Moreover, we make a change of variables
$y_1,y_2$ by $y_1+b_1c_1, y_2+b_2c_2$ respectively.
Then, the system of differential equations for $g$ is given by
$I' \cdot g = 0$, where
\begin{align*}
I' =\langle &A_{11}'=x_{11}\pd{11}-x_{11}\pd{1}^2-a_{11}^3,\;
    A_{12}'=x_{12}\pd{12}-x_{12}\pd{1}\pd{2}-a_{12}^3,\\
& A_{22}'=x_{22}\pd{22}-x_{22}\pd{2}^2-a_{22}^3, \\
& B=\pd{1}^2+\pd{2}^2-r^2, \\  
& C_{12}'=x_{12}\pd{1}^2+2(x_{22}-x_{11})\pd{1}\pd{2}-x_{12}\pd{2}^2
  +(y_2+b_2c_2)\pd{1}-(y_1+b_1c_1)\pd{2}, \\
& E'=r\pd{r}-2(x_{11}\pd{1}^2+x_{12}\pd{1}\pd{2}+x_{22}\pd{2}^2)\\
& \qquad\qquad\qquad\qquad\qquad\qquad-((y_1+b_1c_1)\pd{1}+(y_2+b_2c_2)\pd{2})-1-d^3 \rangle.
\end{align*}
We note that ${\rm rank}(I) = {\rm rank}(I')$ holds for any set of values
of the constants.

Secondly,
we calculate ${\rm in}_{(-w,w)}(I')$ for the weight vector
$w=(0,1,0,0,0,0)$ 
where each weight stands for
the variables $\pd{11},\pd{12},\pd{22},\pd{1},\pd{2},\pd{r}$ respectively.
In order to perform Buchberger's algorithm
with respect to the $(-w,w)$-weight order,
we need to consider the homogenized system $I'^h$ for $I'$:
\begin{align*}
I'^h =\langle &A_{11}'^h=hx_{11}\pd{11}-x_{11}\pd{1}^2-\underline{a_{11}^3},\;
    A_{12}'^h=hx_{12}\pd{12}-x_{12}\pd{1}\pd{2}-\underline{a_{12}^3},\\
& A_{22}'^h=hx_{22}\pd{22}-x_{22}\pd{2}^2-\underline{a_{22}^3}, \\
& B=\pd{1}^2+\pd{2}^2-\underline{r^2}, \\  
& C_{12}'^h=x_{12}\pd{1}^2+2(x_{22}-x_{11})\pd{1}\pd{2}-x_{12}\pd{2}^2
  +(y_2h+b_2c_2)\pd{1}-(y_1h+\underline{b_1c_1)\pd{2}}, \\
& E'^h=hr\pd{r}-2(x_{11}\pd{1}^2+x_{12}\pd{1}\pd{2}+x_{22}\pd{2}^2)\\
& \qquad\qquad\qquad\qquad\qquad-((hy_1+b_1c_1)\pd{1}+(hy_2+b_2c_2)\pd{2})-h^3-\underline{d^3} \rangle.
\end{align*}
We denote by $<_{(-w,w,0)}^h$
an order in the homogenized Weyl algebra
which compares the total degree first, $(-w,w,0)$-degree second,
otherwise we apply the following block order as a tie breaker:
$d \gg r \gg \{a_{11}, a_{12}, a_{22}\}
\gg \{b_1 > b_2\} \gg \{c_1,c_2\} \gg \{y_1,y_2\}
\gg \pd{r} \gg \{\pd{12}\} \gg \{\pd{11}, \pd{22} \} \gg \{\pd{1},\pd{2}\}
\gg \{x_{12}\} \gg \{x_{11},x_{22}\} \gg h$.
Here, the symbol $>$ represents the lexicographic order.
The underlined parts in $I'^h$ are initial terms
with respect to $<_{(-w,w,0)}^h$.
They are pairwise coprime and their commutators
are equal to zero except $[B,E'^h] = -2hB$.
From Lemma \ref{lem:nc-cri-d},
we conclude that the set $\{A_{11}'^h,A_{12}'^h,A_{22}'^h,B,C_{12}'^h,E'^h\}$
is a Gr\"obner basis of $I'^h$ in
$D^h[a_{11},a_{12},a_{22},b_1,b_2,c_1,c_2,d]$
with respect to $<_{(-w,w,0)}^h$.
In other words,
the transformation of dependent and independent variables
gives us the Gr\"obner basis without
adding new elements.
Dehomogenizing $I'^h$,
we obtain
\begin{align*}
{\rm in}_{(-w,w)}(I') =
\langle &\tilde{A}_{11}'=x_{11}\pd{11}-x_{11}\pd{1}^2-a_{11}^3,\;
    \tilde{A}_{12}'=x_{12}\pd{12}-a_{12}^3,\\
& \tilde{A}_{22}'=x_{22}\pd{22}-x_{22}\pd{2}^2-a_{22}^3, \\
& B=\pd{1}^2+\pd{2}^2-r^2, \\  
& \tilde{C}_{12}'=2(x_{22}-x_{11})\pd{1}\pd{2}
  +(y_2+b_2c_2)\pd{1}-(y_1+b_1c_1)\pd{2}, \\
& \tilde{E}'=r\pd{r}-2(x_{11}\pd{1}^2+x_{22}\pd{2}^2)\\
& \qquad\qquad\qquad\qquad-((y_1+b_1c_1)\pd{1}+(y_2+b_2c_2)\pd{2})-1-d^3 \rangle.
\end{align*}
In this calculation, we regard $a_{11},a_{12},a_{22},b_1,b_2,c_1,c_2,d$
as ring variables with the weight 0.
As in the proof of Theorem \ref{th:grdeformation},
the equation above
holds when $a_{11}, \ldots, d$ are specialized to generic complex numbers.
The holonomic rank of the $(-w,w)$-initial ideal $\tilde{I}':={\rm in}_{(-w,w)}(I')$
coincides with that of the diagonal system transformed
by the same change of variables for $I'$ from $I$.
The holonomic rank of $\tilde{I}'$ agrees with
that of
\begin{align*}
\tilde{I} =
\langle &\tilde{A}_{11}=\underline{\pd{11}}-\pd{1}^2,\;
    \tilde{A}_{12}=\underline{\pd{12}},\\
& \tilde{A}_{22}=\underline{\pd{22}}-\pd{2}^2, \\
& B=\underline{\pd{1}^2}+\pd{2}^2-r^2, \\  
& \tilde{C}_{12}=2(x_{22}-x_{11})\underline{\pd{1}\pd{2}}
  +y_2\pd{1}-y_1\pd{2}, \\
& \tilde{E}=r\underline{\pd{r}}-2(x_{11}\pd{1}^2+x_{22}\pd{2}^2)
-(y_1\pd{1}+y_2\pd{2})-1 \rangle.
\end{align*}
Hence, we obtain the following inequality:
$${\rm rank}(I) = {\rm rank}(I')
\geq {\rm rank}({\rm in}_{(-w,w)}(I'))=
{\rm rank}(\tilde{I}') = {\rm rank}(\tilde{I}).$$

Finally, we show that ${\rm rank}(\tilde{I}) = 4$.
Proposition \ref{pr:fbgr2} tells us that the set
\begin{align*}
\{ &\tilde{A}_{11}, \tilde{A}_{12}, \tilde{A}_{22},
B, \tilde{C}_{12}, \tilde{E},\\
&\left. D_2 = 2(x_{22}-x_{11}) \left(\underline{\pd{2}^3}-r^2\pd{2}
 -\frac{y_2\pd{1}^2-y_1\pd{1}\pd{2}-\pd{2}}{2(x_{22}-x_{11})} \right) \right\}.
\end{align*}
is a Gr\"obner basis of $R\tilde{I}$ with respect to
the block order $\{\pd{r}\} \gg \{\pd{11}>\pd{22}\} \gg \{\pd{1}>\pd{2}\}$,
where the tie breaker $>$ represents the graded lexicographic order.
The set of standard monomials is $\{1,\pd{1},\pd{2},\pd{2}^2 \}$.
It means that the holonomic rank of $\tilde{I}$ is $4$.

\end{document}